# S. Barry Cooper: Incomputability After Alan Turing

## 1. Living in a computable world

Those of us old enough may remember being fascinated by George Gamow's popular books on mathematics and science - with the most famous being "One Two Three ... Infinity" . In it, he got us to imagine living on the surface of a 2-dimensional balloon, with only 2-dimensional experience of the surface. And then got us to understand how we might detect its 3-dimensional curved character via purely 2-dimensional observations. Here is his picture from page 103 of the 1961 edition:

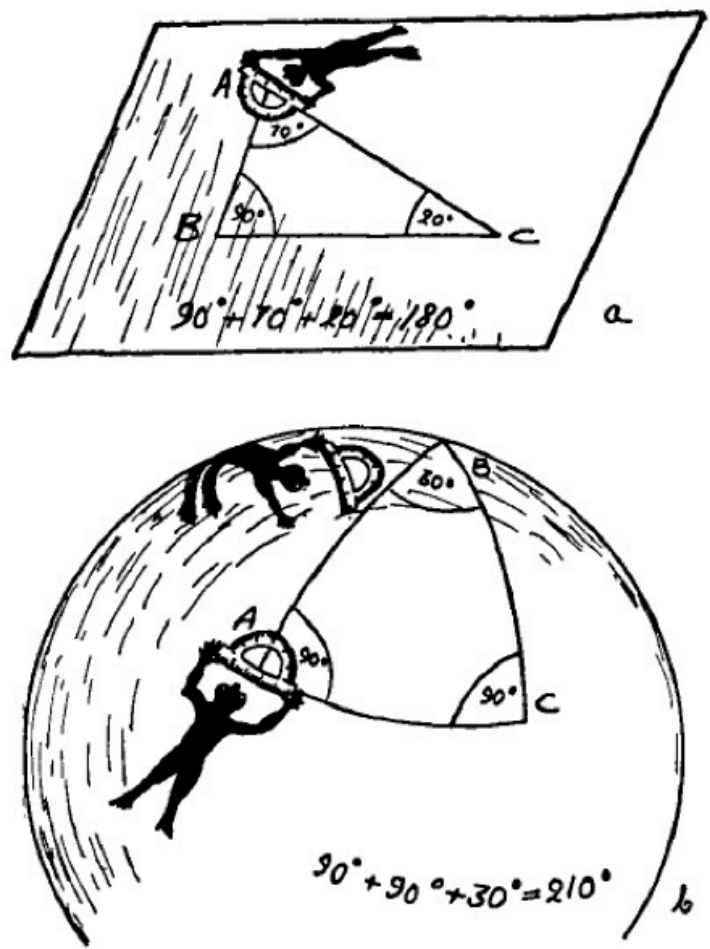

Algorithms, as a way of traversing our four dimensions have been with us for literally thousands of years. They provide recipes for the control and understanding of every aspect of everyday life. Nowadays, they appear as computer programs. Algorithms, or computer programs, can be thought of as a kind of causal dimension all of their own. And then the questions arise: Is there a causal dimension that is not algorithmic? And does it *matter* if there is?

Notice that George Gamow's example showed that on the one hand it was tricky to live in two dimensions and find evidence of a third. And that it *did* matter, precisely because we *could* find that evidence. Of course, if we took the mathematical model presented by the picture, the missing dimension becomes clear to us - we have an overview. But observe, while the mathematical overview gives us a better understanding of the nature of curved space, it does not tell us that the model is relevant to our world. We *still* need to look at the triangle from within the two dimensional world to match up reality and mathematics, and be able to apply the full power of the model.

Back in the 1930s, people like Kurt Gödel, Stephen Kleene, Alonzo Church and Alan Turing did build mathematical models of the computable dimension of causal relations. And this enabled Church and Turing to get outside it, and by examining the model identify the new dimension of incomputability.

A specially important part of what the 24 year-old Alan Turing did was to base his investigation of the extent of the computable on a new *machine-like* model. And the Turing machine was to make him famous in a way he could never have foreseen. He had the idea of using Gödel's coding trick to turn his Turing machine programs into data the machine could compute from - and hence was born the *universal* Turing machine, able to take a code for *any* other given machine as input data, and to compute exactly like it. The universal machine *stored* programs: and so gave us an understanding of the modern stored-program computer before anyone had even built a real one.

This caused all sorts of problems. Just like with George Gamow's example, it was easy to get the mathematical overview. The problem was to match it up with reality. And this was a problem with a practical aspect. Even a toy avatar of the abstract machine was hard to make. The engineers eventually came up with clever solutions to the problem, the EDVAC in Pennsylvania, the Manchester 'Baby', Maurice Wilkes' EDSAC, and the Pilot ACE growing out of Turing's own attempts to build a computer at the NPL. But to this day, there are engineers who don't forgive (or understand even) Turing's reputation as the 'inventor' of the computer.

More importantly, the computer changed every aspect of our lives, and strengthened our experience of living in a computable world. Incomputability became a mathematical oddity, a playground for researchers who were not too concerned about real world significance, but liked doing hard mathematics with a distinct feel of reality. Of *course* it felt like reality. It was real numbers connected by computable relations - a bit like a well-behaved scientific world of information, structured by computable causal relations.

## 2. A short history of incomputability

Computability has always been with us. The universe is full of it - natural laws whose computability enables us to survive in the world; animal and human behaviour guided by biological and learned algorithms; computable natural constants such as $\pi$ and e. The algorithmic content gives the mathematics of nature its infinitary character, and opens the door to incomputability. Richard Feynman may have said (Feynman 1982) "It is really true, somehow,
that the physical world is representable in a discretized way, and ... we are going to have to change the laws of physics": but real numbers persist in the mathematics of the real world.

There have ever been doubts about our ability to make sense of causality. And adding deities may emphasise the problem. It certainly goes back to the 11th century and Al-Ghazali's "The Incoherence of the Philosophers", and is traceable through Hume and Berkeley, arriving, for instance, at the modern interest in emergent phenomena. According to the Oxford English Dictionary (1971 edition), the first recorded use of the word 'incomputable' goes back to 1606 - some 40 years even before 'computable'. The term only acquired its precise meaning in the 1930s, with the formulation of a number of different models of what it means for a function over the natural numbers to be computable. And, as we already mentioned, it was these enabled Church and Turing to get their examples of incomputable objects. The key observation, captured in what we now know as the *Church-Turing Thesis*, is that there is a robust intuitive notion of computability to which all our different formalisms converge. It was Turing's carefully argued 1936 paper, based on the Turing machine model, that convinced Gödel it was true.

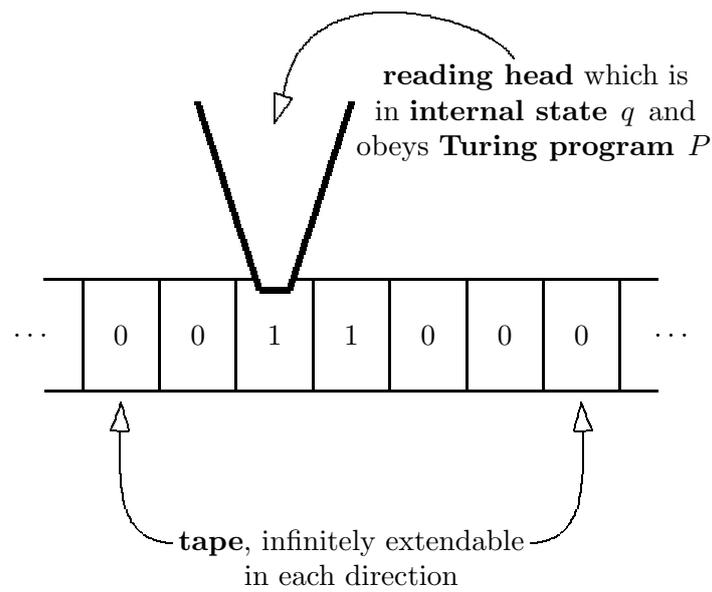

As Gödel's friend Hao Wang recounts (*Reflections on Kurt Gödel*, 1987, p.96):

> Over the years G habitually credited A. M. Turing's paper of 1936 as the definitive work in capturing the intuitive concept [of computability], and did not mention Church or E. Post in this connection. He must have felt that Turing was the only one who gave persuasive arguments to show the adequacy of the precise concept ... In particular, he had probably been aware of the arguments offered by Church for his 'thesis' and decided that they were inadequate. It is clear that G and Turing (1912-1954) had great admiration for each other, ...

Within mathematics, Turing's paper was a blow to David Hilbert's view, famously expressed on September 8th in a Königsberg address, that:

> For the mathematician there is no Ignorabimus, and, in my opinion, not at all for natural science either. . . . The true reason why [no one] has succeeded in finding an unsolvable problem is, in my opinion, that there is no unsolvable problem. In contrast to the foolish Ignorabimus, our credo avers:
>
> We must know, We shall know.

The universality of Turing's machine meant it had to implement a lot of faulty programs. Some programs would lead to computations that never stopped. And one could not tell in general which inputs led to a proper computation. The Halting Problem for the Universal Turing Machine (UTM) turns out to be unsolvable. Or, if you want the full story: The set of inputs which lead to a terminating computation of the UTM is *computably enumerable* - you can progressively set in motion an array of all possible computations and observe which ones output a result. This enables you to enumerate the inputs on which the machine halts. But - this set is *not computable,* since you can *never be certain* that a computation in progress will not one day succeed.

More dramatically, all sorts of mathematical theories are capable of 'talking about' our UTM. Turing used natural numbers to code up the activities of the machine, using the trick Gödel had earlier used to enable Peano arithmetic to talk about itself. Turing's discovery was that *any* reasonably strong mathematical theory was *undecidable* - that is, had an incomputable set of theorems. In particular, Turing had a proof of what became known as *Church's Theorem* - telling us that there is no computer program for testing a statement in natural language for logical validity. Since then a huge number of undecidable theories have been found.

The drawback to this powerful technique for proving the existence of many natural incomputable sets, is that there are no other good techniques for proving incomputability. And John Myhill showed in 1955 that the known natural examples of incomputable objects tended to be all the same; computationally, just notational translations of each other. This would not have mattered, if it were not for the fact that the so-called Turing universe of incomputable sets turned out to have a very rich and mathematically challenging structure. And if this were to be embodied in the real world to any extent, most of the computational character of its embodiment would be hidden from us. If we managed to solve a problem, all well and good. If we failed, we might never know whether that bit of the world had an incomputable character different from that of the halting problem, or might be computable via some program which eludes us so far.

This simple difficulty with *recognising* mathematical incomputability was explained by results which told us such recognition was itself a highly incomputable problem. And here were all the ingredients for a parting of the ways between the mathematics and real-world concerns. Turing himself made his last great contribution to the logic of computability theory in his amazing 1939 paper, based on his work with Alonzo Church in Princeton. After that his work was more obviously rooted in reality than the earlier overarching abstraction. The mathematicians had the halting problem and its variants, mathematically comprehensive and canonical, a little too grandiose for the everyday ad hoc world. While the endless complications of everyday existence could not be classified - the theory was useless.

For the recursion theoretic period of the mathematics, with its isolation and loss of a sense of real-world remit, see recent papers by Robert I. Soare. Turing's Manchester work on AI, connectionist models, and morphogenesis were pointers to the confusions and inspired anticipations of the shape of things to come.

## 3. Mathematical steps towards an incomputable reality

The 1930s unsolvable problems may be examples of incomputable objects, but their mathematical abstraction appears far from the embodied mathematics of a Newton or Einstein. And the deep and intractable problems from the real world are hard to subject to logical analysis. What makes Turing's work so important is the way it draws out the computability-theoretic core of very different real world mysteries. He had a knack of getting inside structures and imaging their constructive cores in new ways. Typically,

Turing does not *apply* mathematics, he *builds* it within the context he is exploring. When Einstein says (Albert Einstein, p.54, 'Out of My Later Years, 1950):

> When we say that we understand a group of natural phenomena, we mean that we have found a constructive theory which embraces them.

Turing takes it literally. And mathematics for a Turing often comes mentally embodied. The for many repellent abstraction of Turing's 1939 paper is all about how Gödel's incompleteness theorem plays out in practice. We notice that the true sentence unprovable in Gödel's theory for arithmetic is easily described. And that enlarging our theory by adding the sentence gives us a larger theory, which has a similarly described unprovable sentence. For Turing this was the seed for a computably iterated process of enlargement. And using Kleene's computable ordinals, one could extend the process transfinitely. Later, logicians such as Solomon Feferman and Michael Rathjen would enhance Turing's ladder into the incomputable, and proof-theoretically reach levels far beyond where Turing had got to. But in a sense, Turing had found out what he wanted to know. Tucked away amongst the mountains of abstraction is his characteristically candid take on what he had done (Turing 1939, pp.134–5):

> Mathematical reasoning may be regarded ... as the exercise of a combination of ... intuition and ingenuity. ... In pre-Gödel times it was thought by some that all the intuitive judgements of mathematics could be replaced by a finite number of ... rules. The necessity for intuition would then be entirely eliminated. In our discussions, however, we have gone to the opposite extreme and eliminated not intuition but ingenuity, and this in spite of the fact that our aim has been in much the same direction.

And mathematically, we have a brilliant analysis of how we may constructively navigate our way through a phase transition, but storing up a level of incomputability arising from the ad hoc nature of the route. There is a degree of arbitrariness in a climber's choice of hand and foot holds, which gives the climb more than algorithmic interest. This ties in with an irreversibility of computation noted by people such as Prigogine in quite different contexts.

There is a special interest for the mathematician in Turing's analysis of an incomputable route to a computable outcomes. It fits nicely with our experience of finding our theorems creatively, and then digging out an algorithmic route to them: something with a memetic character which circulates the community like a virus. A basic rule for lifting small truths to bigger ones, a logical counterpart of the more visceral causality of nature, is mathematical induction. For proof theorists induction plays a key role. They categorise theorems according to the level of complexity of the induction used in the proof. Most theorems turn out to be proof-theoretically very simple. Now that we have a proof of Fermat's Last Theorem, the logician Angus MacIntyre has been able to outline it within first-order arithmetic. This view of the proof involves simple incremental accretions of truth. The discovery of the proof was something very different. And so is our understanding of it.

We are beginning to see a pattern - literally: simple rules, unbounded iteration, emergent forms - defined at the edge of computability. This is just what Turing later observed in nature and mathematically tried to capture.

Another hugely important mathematical tool to take on our explorations of the incomputable is the oracle Turing machine, also tucked away, on one page of Turing's

1939 paper. The idea was to allow the machine to compute relative to a given real, which may or may not be computable. If one looked at the function computed using this oracle, you could see it as a computation of the function from the given oracle real. And it is a small step to summarise what the machine does as computing one real from another. It delivers us a computational model within which to fit basic computable laws of nature - namely, most of what underlies our knowledge of how the world works. Harking back to the computable numbers of 1936, the oracle machine computes a real number, but *from another* real number. We can allow the machine to compute relative to different oracles, then acknowledging the higher type nature of the computational process by calling the mappings computed by oracle machines *Turing functionals*. These, acting over the reals, give us the *Turing universe*.

In fact, Turing was not notably interested in the mathematical development of his oracle machines, despite his preoccupation with computers that interacted. It was left to another seminal figure, Emil Post, to gather together equivalence classes of reals - or *degrees of unsolvability* - which were computable from each other. And then, using an ordering induced by the ordering of reals via Turing functionals, obtain a structure which has become known simply as the *Turing degrees*.

The three things you need to know about this structure are: Firstly that it is very complex; Secondly, that if you take some scientific domain described in terms of real numbers and computable laws over them, then it is embeddable in the Turing universe - so that the structure of the corresponding restriction of the Turing universe tells you something about the causal structure of the real world (the causal third dimension of the Gamow-like 2-dimensional person we met earlier); And thirdly, we can view this model as a terrain in which computation can be hosted, but in which information takes a leading role, structuring the form of the computation - a *re-embodiment* of computation.

## 4. Messages from the real world

Turing's final years in Manchester saw both personal eclipse, and the sowing of the seeds of a later renaissance - a growing renown and increasing scientific impact. The fame is a mixed blessing. All those webpages and articles in periodicals with misleading information about what Turing did and did not do. On the other hand, having mathematicians like Turing and Gödel in the public eye are good for basic science. As far as the science itself goes, Turing's work has been powerfully influential in a piecemeal way, with a number of different fields laying claim to particular bits of Turing. In the October 2004 Notices of the AMS, Lenore Blum wrote a nice description of the dichotomy between "two traditions of computation" (in *Computing over the Reals: Where Turing Meets Newton*, pp. 1024-1034):

> The two major traditions of the theory of computation have, for the most part, run a parallel non-intersecting course. On the one hand, we have numerical analysis and scientific computation; on the other hand, we have the tradition of computation theory arising from logic and computer science.

Turing's 1948 Rounding-off errors in matrix processes was influential in the former, and the 1936 Turing machine paper in the latter.

But there is a growing appreciation of the coherence of approach represented by these different contributions. On the one hand we have a computational world over which we have control; in the other we must live with approximations and errors. As we move up the informational type structure from discrete to continuous data, we lose the sure footholds, but identify emergent controls at higher levels. Turing had earlier introduced Bayesian codebreaking methods at Bletchley Park, in another adjustment to the realities of extracting form from a complex world.

Turings late great contributions traced the computational content of phase transitions in the real world from two different vantage points. Turing had been interested in the emergence of form in nature from his school days - see his mother's sketch of Alan "Watching the Daisies Grow". A short piece by Peter Saunders in the forthcoming *Alan Turing - His Work and Impact* (edited by Cooper and van Leeuwen) discusses Turing's motivation and background reading in getting interested in morphogenesis. He writes:

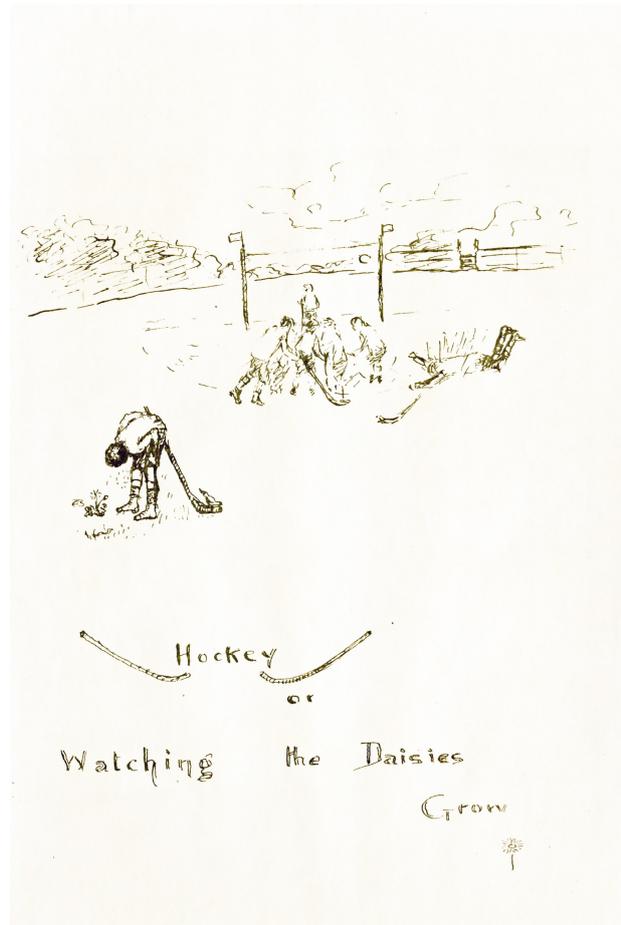

> The obvious question to ask about "The Chemical Basis of Morphogenesis" is why Turing took up the problem at all. Pattern formation, interesting though it may be to biologists, does not look like the sort of fundamental problem that Turing characteristically chose to devote his time and effort to. The answer is simply that he saw it not as a mere puzzle but as a way of addressing what he considered to be a crucial issue in biology. As he said to his student Robin Gandy, his aim was to "defeat the argument from design".

It seems Turing wanted to add a bit more computational convergence to Darwin's theory, and not leave an opportunity for God to tidy up. In the same volume, Philip Maini outlines the way in which Turing produced complex outcomes from simple algorithmic ingredients:

> Alan Turing's paper, "The chemical basis of morphogenesis" (Turing, 1952) has been hugely influential in a number of areas. In this paper, Turing proposed that biological pattern formation arises in response to a chemical pre-pattern which, in turn, is set up by a process which is now known as *diffusion-driven instability*. The genius of this work was that he considered a system which was stable in the absence of diffusion and then showed that the addition of diffusion, which is naturally stabilizing, actually caused an instability. Thus it was the integration of the parts that was as crucial to the understanding of embryological development as the parts themselves - patterns emerged or self-organized as a result of the individual parts interacting. To see how far ahead of his time he was, one has to note that it is only now in the post-genomic era of systems biology that the majority of the scientific community has arrived at the conclusion he came to some 60 years ago.

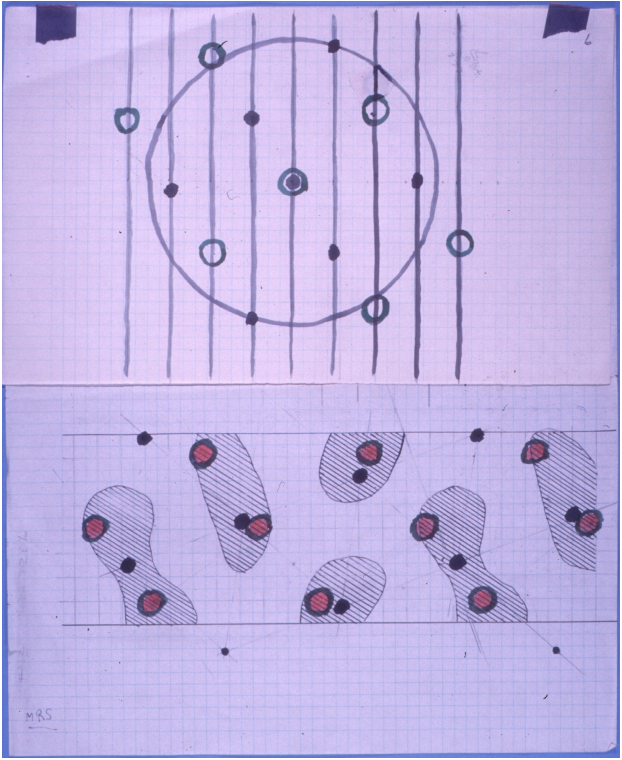

Things have moved on since 1952, but the basic approach retains a powerful influence on the field, hardly noticed by logicians and computer scientists as a group.

More importantly for us, Turing's examples point to general principles underlying emergent phenomena. Turing would have been a teenager at the time the 'British emergentists' such as C. D. Broad, Samuel Alexander, and C. Lloyd Morgan were at their height. Broad and Turing overlapped at Cambridge. As a group they were very prescient in their seeking out examples of the complex arising from simple rules in highly connected environments. Unfortunately, some of their examples of emergence taken from chemistry turned out to be explainable in terms of quantum mechanics. But Turing's differential equations gave us a new sense both of the character and the origins of emergence, and indicated how emergent form night be captured mathematically.

Turing's equations might well have computable solutions. But they pointed to the principle of emergent patterns on animal coats etc. corresponding to definable relations over basic mathematical structure. If one could define a relation in nature, it had a robustness, a tangible presence which one might expect to find observable; just as Turing had brought us to expect observed emergent phenomena to have descriptions. Descriptions which might lead to incomputable sets if at least as complicated as that giving us the halting problem. Somewhere between the tidy abstraction of the universal Turing machine and the mysteries of emergence in nature, we have the fractal family. Well known as analogues of emergent phenomena, they also have their own well-defined mathematics. For various reasons the most informative of these is the Mandelbrot set. As Roger Penrose puts it in his 1994 book on The Emperor's New Mind:

> Now we witnessed ... a certain extraordinarily complicated looking set, namely the Mandelbrot set. Although the rules which provide its definition are surprisingly simple, the set itself exhibits an endless variety of highly elaborate structures.

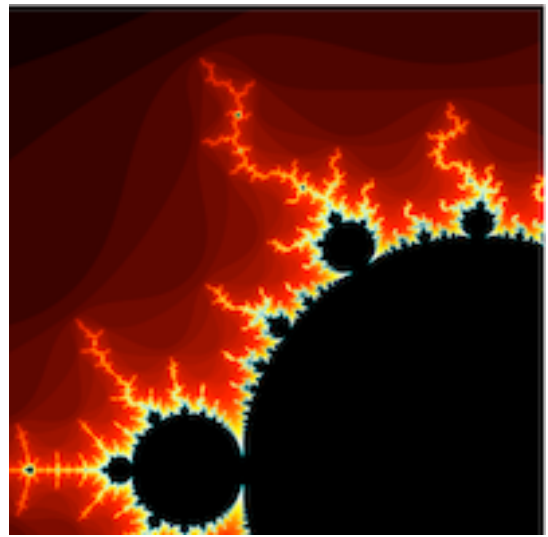

With its definition based on a well-known simple equation over the complex numbers, the quantifier form of its complement can be reduced to something similar to that of the halting set of the UTM. And using this one can simulate it on a computer screen, getting the fascinating range of images we recognise so easily. It has what the the halting problem does not, the visual embodiment of a natural phenomenon. And this enables us to appreciate the higher-order intricacy we encounter as we travel deeper and deeper into this endlessly surprising mathematical object. It is harder to be sure of ones notion of computability in this context. But for the computable analysts, the computability of the Mandelbrot set is still a challenging open problem.

Compared to other levels of the real world that Turing was drawn to, the emergence seems relatively straightforwardly mathematical. One has an objective view of the whole picture, basic rules and emergence of a surprising character - the surprise is one of the criteria for emergence, though there is no proper definition. Mathematically, we tend to look for quantifiers or non linearity in the description based on the basic operations, a kind of association of true emergence with halting problem-like incomputability. At the quantum level, which Turing was always very interested in, but did not live long enough to get a grip on - the problem of pinning down the basic causality is not so easy. We are looking down from above, and are never sure we have the whole picture. And the phase transition from quantum ambiguity to familiar classical world seems to involve not just definability, but a transition from a structure in which there are some prohibitions on simultaneous definition of entities - something familiar to model-theorists - which disappear as one passes from particle physics to other scientific areas. There is experimental evidence that such breakdowns of definability occur in human mentality. The problem we as observers have in this context is not one of viewing from above, but of being trapped inside; although modern neuroscience is augmenting this inner view with a huge amount of useful information.

In his final years, Turing approached brain functionality from two directions: mathematically modelling the physical connectivity of the brain; and via his much better known discussion of intelligent thought. The latter is of special interest for mathematicians, and not just for its relevance to incomputability. Around the same time as Turing was approaching 'intuition and ingenuity' via his hierarchical analysis of the limits of Gödel's theorem, Jacques Hadamard was covering very similar ground from a more sociological viewpoint. Given that the mathematical product is presented algorithmically via a proof, it makes mathematical thinking a good case study for clarifying the 'intuition versus ingenuity' dichotomy. A principal source for Hadamard's 1945 book *The Psychology of Invention in the Mathematical Field* were lectures of Henri Poincare to the Société de Psychologie in Paris in the early part of the 20th century. Here is Hadamard's account of an example of apparent non-algorithmic thinking:

> At first Poincaré attacked [a problem] vainly for a fortnight, attempting to prove there could not be any such function ... [quoting Poincaré]:
> 'Having reached Coutances, we entered an omnibus to go some place or other. At the moment when I put my foot on the step, the idea came to me, without anything in my former thoughts seeming to have paved the way for it ... I did not verify the idea ... I went on with a conversation already commenced, but I felt a perfect certainty. On my return to Caen, for conscience sake, I verified the result at my leisure.'

Many writers focus on the surprise and dissociation from consciously rational thought. What is also striking is the connection with Turing's 1939 paper - the 'perfect certainty' that Poincaré experienced. Did Poincaré have all the details of the proof immediately sorted in his mind? Unlikely. What we understand from Turing's analysis is that there is a process of definition, with a range of codes for proofs emerging. It was one of these codes extracted by  Poincaré on his return on his return to Caen.

As well as the connection back to Turing's 1939 paper, there is the connection forward to his 1952 paper and the final work defining emergence of form in nature. Since 1954, neuroscience is one of the areas in which emergence has taken centre-stage for many researchers. And Turing's connection between definability in terms of both logical structure and physical context is an remarkable anticipation of current thinking.

At the bottom of the process is the basic physical functionality, and Turing made his own groundbreaking contribution to connectionist models of the brain in his unpublished 1948 National Physical Laboratory report on *Intelligent Machinery*. Called unorganised machines by Turing, they were pre-empted by the more famous neural nets of McCulloch and Pitts. See Christof Teuscher's book on *Turing's Connectionism*. These are Teuscher's comments on the history in his commentary on *Intelligent Machinery* in *Alan Turing - His Work and Impact*:

> In his work, Turing makes no reference to McCulloch and Pitts' 1943 paper, nor do they mention Turing's work [on unorganised machines]. It is an open question how much their work influenced each other, yet, we have to assume that they were at least aware of each other's ideas. We hypothesize that both bad timing and the fact that Turing's neurons are simpler and more abstract contributed to his work being largely ignored.

Connectionist models have come a long way since Turing's time. Their physical emulation of the brain does bring dividends. Paul Smolensky, for instance, talks in his 1988 paper *On the proper treatment of connectionism* of a possible challenge to "the strong construal of Church's Thesis as the claim that the class of well-defined computations is exhausted by those of Turing machines."

Of course, it was the celebrated 1950 *Mind* paper on *Computing Machinery and Intelligence* that became one of Turing's three most cited papers. And the Turing Test for intelligence has entered the popular culture, with students walking around with T-shirts proclaiming "I failed the Turing Test".

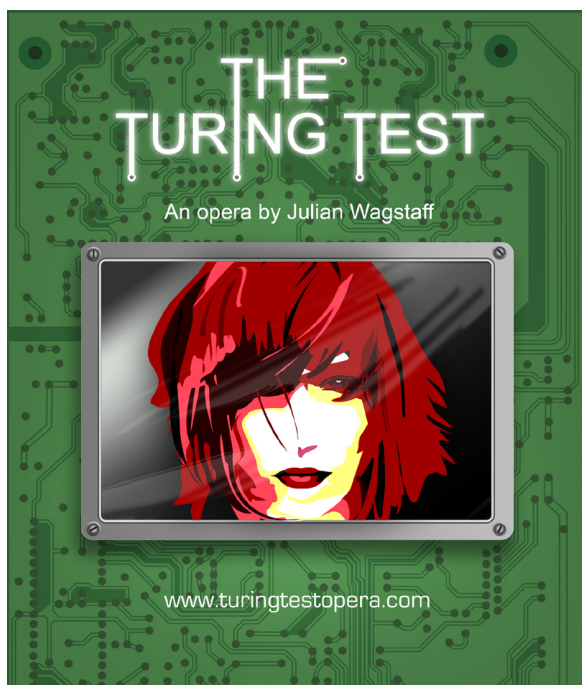

So what is all this to do with the mathematics of incomputability? For all Turing's focus on the logical structure of computation, which has had a huge influence on modern thinking - particularly since Bletchley Park Turing was involved in many ways with *embodied* computation. This involved an implicit recognition of a difference in the way humans and intelligent machines in general, interacted - were *embedded* in - *information*. And this was a necessary attribute of intelligence. People involved with taking artificial intelligence forward have had to take on board this extended, embodied, physical, information respecting, model of the future. We give some brief excerpts from Rodney Brooks' contribution (*The Case For Embodied Intelligence*) to *Alan Turing - His Work and Impact*:

> For me Alan Turing's 1948 paper Intelligent Machinery was more important than his 1950 paper Computing Machinery and Intelligence. ...
>
> For me, the critical, and new, insights in Intelligent Machinery were two fold. First, Turing made the distinction between embodied and disembodied intelligence. ...
>
> Modern researchers are now seriously investigating the emobodied approach to intelligence and have rediscovered the importance of interaction with people as the basis for intelligence. My own work for the last twenty five years has been based on these two ideas.

## 5. Turing points the way past the Turing barrier

Alan Turing's work was incomplete, as it will be for all of us when the final day comes. For Turing it came too early, via an uneasy confluence of algorithm (the UK law of the time) and incomputability, a bizarre piece of history that nobody could have invented for a man who had served mathematics and science - and his country - so well. Some have questioned the description 'computability theory' for the subject Turing co-founded with other luminaries of the period - Godel, Post, Church, Kleene - because it deals primarily with the incomputable.

Turing was a mathematician of his time, who worked from within the world, trying to give mathematical substance to the physical and mental processes. He encountered many interesting things which he left those who came after to take forward. He gave us a basic model of what we understand to be computation. And took observation of it as an organic whole to the next level, discovering incomputability, a definability-theoretic extension of Gödel's incompleteness theorem. He encouraged us to see the universe as something which does compute, and engaged with its features, tracing embodied analogues of his halting problem within biology and neuroscience. He loved the truth, and was open to doubts about his doubts. He thought a machine could not be intelligent if it was expected to be infallible. He saw something different about an embedded computer, and was drawn to the computational mysteries of quantum theory.

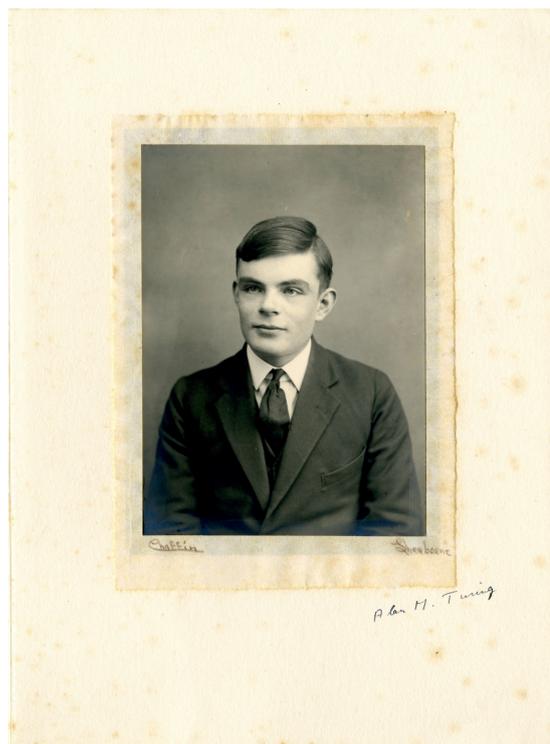

Turing did not live long enough to appreciate Stephen Kleene's investigations of higher-type computability, but would surely have made the connection between the mathematics of incomputability, definability, and the computations that arise from them. He did not see the mathematical theory of randomness flourish, or the yearly award of a prize in his name for a subject he played a founding role in. He never saw the employment of thousands who talked about 'Turing' machines and a 'Turing' Test. The 'incomputable reality' (as Nature recently described it) is still dangerous to inhabit. But we do have a lot to celebrate in 2012.